\documentclass{article}

\usepackage{arxiv}

\usepackage[utf8]{inputenc} % allow utf-8 input
\usepackage[T1]{fontenc}    % use 8-bit T1 fonts
\usepackage{hyperref}       % hyperlinks
\usepackage{url}            % simple URL typesetting
\usepackage{amsfonts}       % blackboard math symbols
\usepackage{lipsum}	
\usepackage{amsmath}
\usepackage{amsthm}
\newtheorem{thm}{Theorem}[section]
 \newtheorem{defini}{Definition}[section]
 \newtheorem{rem}{Remark}[section]
 \newtheorem{lem}{Lemma}[section]

 \numberwithin{equation}{section}

\title{Regularizing effect of absorption terms in singular and degenerate elliptic problems}

\author{
Abdelaaziz Sbai and Youssef El hadfi\\%\thanks{Use footnote for providing further
Laboratory LIPIM\\
National School of Applied Sciences Khouribga\\
Sultan Moulay Slimane University, Morocco\\
  \texttt{yelhadfi@gmail.com}\\
    \texttt{sbaiabdlaaziz@gmail.com}\\
}

\begin{document}
\maketitle
\begin{abstract}
 In this paper we study the existence and regularity of
solutions to the following singular  problem
\begin{equation}
\left\{
\begin{array}{lll}
&-\displaystyle\mbox{div} \big(a(x,u)|\nabla u|^{p-2}|\nabla u|\big) + |u|^{s-1}u =\frac{f}{u^{\gamma}} &\mbox{ in } \Omega \\
&u>0 &\mbox{ in }\Omega \\
&u=0  &\mbox{ on }  \delta\Omega
\end{array}
\right.
%\tag{$P_{a}$}
\end{equation}	
proving that the lower order term  $u|u|^{s-1}$ has some regularizing effects on the solutions in the case of an elliptic  operator with degenerate coercivity.
\end{abstract}

% keywords can be removed
\keywords{Degenerate coercivity, singular non linearity, regularity, entropy solutions, Sobolev spaces.}
%\subclass{35A21, 35B20, 35B25, 35B45, 35D30, 35J60, 35J70, 35J75.}

\section{Introduction}
Let us consider the following problem
\begin{equation}\label{S1}
\left\{
\begin{array}{lll}
&-\displaystyle\mbox{div} \big(a(x,u)|\nabla u|^{p-2}|\nabla u|\big) + |u|^{s-1}u =h(u)f & \mbox{ in } \Omega \\
&u\geq 0 &\mbox{ in }\Omega \\
&u=0  &\mbox{ on }  \delta\Omega
\end{array}
\right.
%\tag{$P_{a}$}
\end{equation}	
where $1<p<N, \,\Omega$ is bounded set in $I\!\!R^{N}$ and $ a:\Omega \times I\!\!R \longrightarrow I\!\!R$	
is a carath\'{e}odory function such that for a.e. $x\in \Omega$ and for every $s\in  I\!\!R$, we have
\begin{equation}\label{S2}
a(x,s)\geq \frac{\alpha}{(1+|s|)^{\theta}}
\end{equation}
\begin{equation}\label{S3}
a(x,s)\leq \beta,
\end{equation}
for some real positive constants $\alpha, \,\beta$ and $\theta.$ Moreover, $f$ is a non negative $L^{m}(\Omega)$ function, with $m\geq 1$ and the term $h:(0,\infty)\longrightarrow (0,\infty)$ is continuous and bounded, such that
\begin{equation}\label{S4}
\exists  c,\, \gamma >0 \mbox{ s.t  } h(s)\leq\frac{c}{s^{\gamma}} \quad \forall s\geq 0,
\end{equation}
for some real number $\gamma$ such that $0\leq \gamma <1.$ Singular problems of this type have been largely studied in the past also for their connection with the theory of non-Newtonian fluids, boundary layer phenomena for viscous fluids and  chemical heterogeneous (see for instance [\cite{14}, \cite{15}]).

Let us briefly recall the mathematical framework concerning problem \eqref{S1} we start with the case  $a(x,u):=a(x), \theta =0$ and  $f$ lies just in $L^{1}(\Omega)$ has been studied in \cite{11}.

Problem \eqref{S1} in the non-singular case $h(u)=1,$ ( problem \eqref{S5})  the author studied the existence and regularity of weak solution to the elliptic problem with degenerate coercivity (see \cite{03}).
\begin{equation}\label{S5}
\left\{
\begin{array}{ll}
-\displaystyle\mbox{ div } \big(a(x,u)|\nabla u|^{p-2}|\nabla u|\big) + |u|^{r-1}u =f, \, in \,\,\,\,\Omega \\
u=0 \,\,\,\, in\,\,\,\,\,  \partial \Omega,
\end{array}
\right.
%\tag{$P_{a}$}
\end{equation}	
in the case where $f \in L^{m}(\Omega)$ with $m\geq 1$ and $\theta >0$. If $p=2$, the problem \eqref{S5} have been treated
in \cite{12}, i.e, in the case of the following problem
\begin{equation}\label{S6}
\left\{
\begin{array}{ll}
-\displaystyle\mbox{div} \big(a(x,u)|\nabla u|\big) + |u|^{r-1}u =f, \, in \,\,\,\,\Omega \\
u=0 \,\,\,\, in\,\,\,\,\,  \partial \Omega,
\end{array}
\right.
%\tag{$P_{a}$}
\end{equation}
the authors studied the  lower order term $|u|^{r-1}u$ in \eqref{S5} and \eqref{S6} that  has the regularizing effects of the solutions in the case where $f\in L^{m}(\Omega)$, with $m\geq 1$ and $\theta \geq 0$.
When p = 2 and the lower-order term does not appear in \eqref{S5}, the existence and regularity of solution to problem \eqref{S5} are proved in \cite{10}. The extension of this work to general case is investigated in \cite{02}.\\
Now we turn our attention recalling some results when the authors had added the singular sourcing term. Problems of p-Laplacien type (i.e $\theta=0$), have been well studied in both the existence and regularity aspects with $f$ having different summbility (see \cite{09}). This frame work has been extended to the problems with a lower order, considering
\begin{equation}\label{S7}
\left\{
\begin{array}{ll}
-\Delta u + u^{s} =\frac{f}{u^{\gamma}}, \, in \,\,\,\,\Omega \\
u>0\,\,\,\, in\,\,\,\,\,   \Omega,
u=0 \,\,\,\, in\,\,\,\,\,  \partial \Omega,
\end{array}
\right.
%\tag{$P_{a}$}
\end{equation}
with $f\in L^{m}(\Omega)$, $m\geq 1$, $0\leq \gamma <1$. Existence and regularity was established in \cite{17}. Recently Olivia \cite{16} have proved the existence and regularity of solution to the problem
\begin{equation}\label{S8}
\left\{
\begin{array}{ll}
-\Delta_{p} u + g(u) =h(u)f, \, in \,\,\,\,\Omega \\
u>0\,\,\,\, in\,\,\,\,\,   \Omega,\\
u=0 \,\,\,\, in\,\,\,\,\,  \partial \Omega,
\end{array}
\right.
%\tag{$P_{a}$}
\end{equation}
$f$ is nonnegative and it belongs to $f\in L^{m}(\Omega)$, $m\geq 1$, for some $0\leq \gamma <1$. While $g(s)$
is continuous, $g(0) = 0 $ and, as $s \rightarrow \infty$, could act as $s^{q}$ with $q\geq -1$ , the p-Laplacian operator is $\Delta_{p} u:=div(|\nabla u|^{p-2}|\nabla u|)$ and  $h$ is continuous, it possibly blows
up at the origin and it is bounded at infinity. \\
In \cite{19}, the authors studied the following degenerate elliptic problem with a singular nonlinearity:
\begin{equation}\label{11}
\left\{
\begin{array}{lll}
&-\displaystyle\mbox{div}( a(x,u,\nabla u)) =fh(u)& \mbox{ in } \Omega \\
&u>0 &\mbox{ in }\Omega \\
&u=0  &\mbox{ on }  \delta\Omega,
\end{array}
\right.
\end{equation}
$f$ is nonnegative and it belongs to $f\in L^{m}(\Omega)$, $m\geq 1$, and $h$ satisfied
the condition in \eqref{S4}. Following this way in this paper, we are interested again in the regularity results. By adding the singular term to the right of \eqref{S5}, we investigate the regularity of solutions of problems of kind \eqref{S1} in light of the influence of some lower order terms.
\par In the study of problem \eqref{S1}, there one to two difficulties, the first one is the fact that, due to hypothesis \eqref{S2} the differential operator $A(u)=div \big(a(x,u)|\nabla u|^{p-2}|\nabla u|\big)$ is not coercive on $W_{0}^{1,p}(\Omega)$, when $u$ is large  (see \cite{13}). Due to the lack of coercivity,  the classical theory for elliptic operators acting between spaces in duality (see  \cite{07}  ) cannot be applied. The second difficulty comes from the right-hand side is singular in the variable $u$. We overcome
these difficulties by replacing operator $A$ by another one defined by means of truncations, and approximating the singular term by non singular one. We will prove in section\eqref{sec3} that these problems admit a bounded $W_{0}^{1,p}(\Omega)$ solution $u_{n}$, $n\in I\!\!N$ by using Schauder's fixed point theorem. In section\ref{sec4} we will get some a priori estimates and convergence results on the sequence of approximating solutions. In the end, we pass to the limit in the approximate problems.

\textbf{Notations:}
In the entire paper $\Omega$ is an open and bounded subset of $ I\!\!R^{N}$, with $N \geq 1$, we denote by $\partial A$ the boundary and by
$|A|$ the Lebesgue measure of a subset $A$ of
$ I\!\!R^{N}$.\\
For any $q>1$, $q^{\prime}=\frac{q}{q-1}$ is the H\"{o}lder conjugate exponent of $q$, while for any $1\leq p<N$, $p^{\ast}=\frac{Np}{N-p}$ is the Sobolev conjugate exponent of $p.$ For fixed $k>0$ we will use of the truncation $T_{K}$ defined as $T_{k}(s)=max \big(-k,min(k,s)\big),$ we will also use the following functions
\begin{equation}\label{N0}
V_{\delta,k}(s)=\left\{
\begin{array}{ll}
1 \,\,\,\,\,\,\,\,      s\leq k \\
\frac{k+\delta-s}{\delta} \,\,\,\,\,\,\,\,\, k <s<k+\delta,\\
0 \,\,\,\,\,\,\,\,\,  s\geq k+\delta,
\end{array}
\right.
\end{equation}
and
\begin{equation}\label{N1}
S_{\delta,k}(s):=1-V_{\delta,k}(s).
\end{equation}
For the sake of implicity we will often use the simplified notation
$$\int_{\Omega}f:=\int_{\Omega}f(x)dx,$$
when referring to integrals when no ambiguity on the variable of integration is possible. If no otherwise specified, we will denote by $c$ serval constants whose value may change from line to line and, sometimes, on the same line. These values will only depend on the data (for instance c can depend on $\Omega,\, \gamma, \,N,\, k,...$) but the will never depend on the indexes of the sequences we will often introduce.
\section{Statement of definitions and the main results }
\subsection{Statement of definitions}
In this context we deal with some class of solutions
\begin{defini}\label{d1}
	A positive function $u$ in $  W_{0}^{1,p}(\Omega)$ is weak solution \eqref{S1} if $h(u)f \in L_{loc}^{1}(\Omega),$ $|u|^{s}$ $\in L^{1}(\Omega) $ and if
	\begin{equation}\label{S5'}
	\int_{\Omega} a(x,u)|\nabla u|^{p-2}\nabla u\nabla \varphi dx + \int_{\Omega}|u|^{s-1}u \varphi=\int_{\Omega}fh(u)\varphi\,\,\,\,\forall \varphi \in C_{c}^{1}(\Omega).
	\end{equation}
\end{defini}
\begin{defini}
	A measurable function $u$ is an entropy solution to problem \eqref{S1} if  $ |u|^{s}\in L^{1}(\Omega), \, h(u)f \in L_{loc}^{1}(\Omega), \,T_{k}(u)\in W_{0}^{1,p}(\Omega)$ for every $k>0$ and
	\begin{equation}\label{T1}
	\int_{\Omega} a(x,u)|\nabla u|^{p-2}\nabla u\nabla T_{k}(u-\varphi)dx+\int_{\Omega}|u|^{s-1}u T_{k}(u-\varphi)dx\leq \int_{\Omega}fh(u) T_{k}(u-\varphi)dx,
	\end{equation}
	for every $\varphi\in W_{0}^{1,p}(\Omega)\cap L^{\infty}(\Omega).$	
\end{defini}	
\begin{defini} The Marcinkwicz space $M^{s}(\Omega)$, $s>0$, consists of all measurable functions $v:\Omega \longrightarrow  I\!\!R$ that satisfy the following condition: there  exists $c>0$ such that
	$$ meas\{|v|\geq k\} \leq \frac{c}{k^{s}} \,\,\,\,\mbox{for all}\,\,\,k>0.$$
	If $|\Omega|<\infty$ and $ 0<\epsilon<s-1$, we can show that $L^{s}(\Omega) \subset M^{s}(\Omega)\subset L^{s-\epsilon}(\Omega).$
	%$$p_{0}:=1+\frac{(1+\theta-\gamma)(N-1)}{N}.$$
	\end{defini}
Let
\begin{equation}\label{S6'}
p_{0}:=1+\frac{(1+\theta-\gamma)(N-1)}{N}
\end{equation}
\subsection{ Statement of  the main results }
The main results of this paper are stated as follows:
\begin{thm}\label{thm1}
	Let $f\in L^{m}(\Omega)$, $m>1$, $1<p<N$. Then
	\begin{itemize}	
		\item[i)] If $s\geq \frac{1+\theta- \gamma}{m-1}$, then there exists a distributional solution $u$ to problem \eqref{S1} such that
		$$ u\in  W_{0}^{1,p}(\Omega)\cap L^{ms+\gamma}(\Omega).$$	
		\item[ii)]  If $\frac{1+\theta -\gamma}{pm-1}<s< \frac{1+\theta- \gamma}{m-1}$, then there exists a distributional solution $u$ to problem  \eqref{S1} such that\\ $$u^{ms+\gamma}\in L^{1}(\Omega)\,\,\,\,\mbox{and}\,\,\,\,u\in W_{0}^{1,\sigma}(\Omega)\,\,\,\,, 1< \sigma=\frac{pms}{1+\theta+s-\gamma}. $$
		\item[iii)]If $0< s\leq \frac{1+\theta -\gamma}{pm-1}$, then there exists an entropy solution u to problem \eqref{S1} such that
		$$u^{ms+\gamma}\in  L^{1}(\Omega)\,\,\,\,\,\mbox{and}\,\,\,\,| \nabla u|\in M^{\frac{pms}{1+\theta+s-\gamma}}(\Omega).$$
	\end{itemize}	
\end{thm}
\begin{rem}
	If $p=2$ and $\gamma=0;$ the result of Theorem \ref{thm1} coincides with regularity results in the case of an elliptic operator with degenerate coercivity ( see \cite{12}, Theorem 1.5).
\end{rem}
\begin{thm}\label{t2}
	Let $f\in L^{m}(\Omega), \,m>1,\, p_{0}<p<N.$ Then
	\begin{itemize}	
		\item[i)] If $0<s\leq \frac{N(1-\gamma )+\gamma }{m(N-1)}$, then there exists a distributional solution $u$ to problem  \eqref{S1}   such that\\   $u^{ms+\gamma}\in L^{1}(\Omega)$ and $u\in W_{0}^{1,\sigma}(\Omega)$, where $1< \sigma=\frac{N[p+s(m-1)-1-\theta+\gamma]}{N+s(m-1)-1-\theta+\gamma}. $\\
		\item[ii)]If $ s\geq \frac{N(1-\gamma )+\gamma }{m(N-1)}$, then item(ii) of Theorem \ref{thm1} holds.
	\end{itemize}	
\end{thm}
\begin{rem}
	If $\gamma=0$; the result of Theorem \ref{thm1} coincides with regularity results in the case of an elliptic operator with degenerate coercivity ( see \cite{03}, Theorem 3)
	and Theorem \ref{t2}  coincides with (\cite{03}, Theorem 4).
\end{rem}
\begin{thm}\label{t3}
	Let $f\in L^{1}(\Omega),$ $1<p<N.$ Then
	\begin{itemize}	
		\item[a)] If $s>\frac{1+\theta -\gamma}{p-1},$ then there exists a distributional solution $u$ to problem \eqref{S1} such that \\  $u^{s+\gamma}\in L^{1}(\Omega)$ and $u\in W_{0}^{1,r}(\Omega)\cap L^{s+\gamma}(\Omega)$,where $1< r<\frac{ps}{s+1+\theta-\gamma}.$
		\item[b)]If $ 0<s\leq \frac{1+\theta -\gamma}{p-1},$ then there exists an entropy solution u to problem \eqref{S1} such that $$u^{s+\gamma}\in L^{1}(\Omega)\,\,\,\,\mbox{and}\,\,\,\,|\nabla u|\in M^{\frac{ps}{s+1+\theta-\gamma}}.$$
	\end{itemize}	
\end{thm}
\begin{rem}
	If $p=2$ and $\gamma=0$; the result of Theorem \ref{t3} coincides with regularity results in the case of an elliptic operator with degenerate coercivity ( see \cite{12}, Theorem 1.4).
\end{rem}
\begin{thm}\label{t4}
	Let $f\in L^{1}(\Omega),$ $p_{0}<p<N.$ Then
	\begin{itemize}	
		\item[1)] If $0<s\leq \frac{N(1-\gamma)+\gamma}{N-1}$, then there exists a distributional solution $u$ to problem \eqref{S1} such that \\   $u\in W_{0}^{1,r}(\Omega)$, where $1< r<\frac{N[p-\theta-1+\gamma]}{N-\theta+\gamma-1}.$\\
		\item[2)]If $ \frac{N(1-\gamma)+\gamma}{N-1}<s< \frac{N(p-1-\theta )+p\gamma}{N-p}$, then item (b) of theorem \ref{t3} holds.\\
		\item[3)] If $ s\geq  \frac{N(p-1-\theta )+p\gamma}{N-p}$, then item (a) of Theorem \ref{t3} holds.
		\end{itemize}	
\end{thm}
\begin{rem}
	If $\gamma=0;$ the result of Theorem \ref{t3} coincides with regularity results in the case of an elliptic operator with degenerate coercivity ( see \cite{03}, Theorem 1) and Theorem \ref{t4}  coincides with (\cite{03}, Theorem 2).
\end{rem}		
\begin{rem}
	In Theorem$\ref{t4},$ we have $\frac{N(1-\gamma)+\gamma}{N-1}<s< \frac{N(p-1-\theta)+p\gamma}{N-p}$ is meaningful, because
	$p>p_{0},$ besides $s\leq \frac{N(1-\gamma)+\gamma}{N-1} \Rightarrow$ $s+\gamma\leq r^{\ast}$, for $r\geq 1.$ In Theorem \ref{t2},
	$s\leq \frac{N(1-\gamma)+\gamma}{m(N-1)}\Rightarrow $ $ms+\gamma\leq \sigma^{\ast},$ for $\sigma \geq 1.$
\end{rem}
\section{A priori estimates and Preliminary facts}\label{sec3}
Let us introduce the following scheme of approximation
\begin{equation}\label{A1}
\left\{
\begin{array}{ll}
-\displaystyle\mbox{div} \big(a(x,T_{n}(u_{n}))|\nabla u_{n}|^{p-2}|\nabla u_{n}|\big) + |u_{n}|^{s-1}u_{n}=h_{n}(u_{n})f_{n}, \, in \,\,\,\,\Omega \\
u_{n}=0 \,\,\,\, in\,\,\,\,\,  \partial \Omega,
\end{array}
\right.
%\tag{$P_{a}$}
\end{equation}
where $f_{n}=T_{n}(f).$ Moreover, defining  $h(0):= \lim_{s\rightarrow0}h(s),$ we set
\begin{equation}\label{A2}
h_{n}(s)=\left\{
\begin{array}{ll}
T_{n}(h(s)) \,\,\,\, \mbox{for} \,\,\,\,\,\, s>0,\\
\min(n,h(0)) \,\,\,\,\,\,\,\,  \mbox{otherwise}.
\end{array}
\right.
%\tag{$P_{a}$}
\end{equation}
The right hand side of \eqref{A1} is non negative, that $u_{n}$ is non negative. The existence of weak solution $u_{n}\in W_{0}^{1,p}(\Omega)$ is guaranteed by the following lemma.
\begin{lem} \label{lem1}
	Problem \eqref{A1} has a non negative solution $u_{n}$ in $  W_{0}^{1,p}(\Omega),$ such that
	\begin{equation}\label{l01}	
	\int_{\Omega} |u_{n}|^{ms+\gamma} dx \leq c\int_{\Omega} |f|^{m} dx
	\end{equation}
	and the solution $u_{n}$ satisfies
	\begin{equation}\label{l02}
	\int_{\Omega} a(x,T_{n}(u_{n}))|\nabla u_{n}|^{p-2}|\nabla u_{n}|\nabla \varphi dx + \int_{\Omega}|u_{n}|^{s-1}u_{n} \varphi=  \int_{\Omega}f_{n}h_{n}(u_{n})\varphi,
	\end{equation}
	where $0\leq \gamma <1$ and $\varphi$ in $W_{0}^{1,p}(\Omega)\cap L^{\infty}(\Omega).$
\end{lem}
\textbf{Proof}.
This proofs derived from Schauder's fixed point argument in \cite{01}.
For fixed  $n \in I\!\!N$ let us define a map
$$ G:L^{p}(\Omega) \rightarrow L^{p}(\Omega),$$
such that, for any  $v$ be a function in $L^{p}(\Omega)$ gives the weak solution $w$ to the following problem
\begin{equation}\label{l1}
-div (a(x,T_{n}(w))|\nabla w|^{p-2}|\nabla w|) + |w|^{s-1}w =  f_{n}h_{n}(v).
\end{equation}
The existence of a unique $ w\in  W_{0}^{1,p}(\Omega)$ corresponding to a $v\in L^{p}(\Omega)$ follows from the classical  result of [\cite{02}, \cite{07}].
Moreover, since the datum $f_{n}h_{n}(v)$ bounded, we have that $w\in L^{\infty}(\Omega)$ and there exists a positive constant $d_{1},$ independents of $v$ and $w$ (but possibly depending in $n$), such that\\ $||w||_{L^{\infty}(\Omega)}\leq d_{1}.$ Again, thanks to the regularity of the datum $f_{n}h_{n}(v),$ we have can choose $w$ as test function in the weak formulation \eqref{31}, we have
\begin{equation}\label{l2}
\displaystyle\int_{\Omega} a(x,T_{n}(w))|\nabla w|^{p-2} \nabla w \nabla w + \displaystyle\int_{\Omega}|w|^{s-1}w.w  = \displaystyle\int_{\Omega}  f_{n}h_{n}(v)w,
\end{equation}
then, it follows from \eqref{S2}
$$
\alpha\displaystyle\int_{\Omega} \frac{|\nabla w |^{p}}{(1+n)^{\theta }}dx \leq n^{2}\int_{\Omega} |w| dx,
$$
using the Poincar\'{e} inequality we have
$$
\displaystyle\int_{\Omega} \frac{|\nabla w |^{p}}{(1+n)^{\theta }}dx \leq \frac{c_{1}}{\alpha} n^{2}\int_{\Omega} |\nabla w| dx,
$$
then
\begin{equation}
\displaystyle\int_{\Omega} |\nabla w |^{p} dx \leq  \frac{c_{1}}{\alpha}(1+n)^{\theta+2}\displaystyle\int_{\Omega}  |\nabla w | dx
\leq c(n,\alpha) |\Omega |^{\frac{1}{p'} } \left(\displaystyle\int_{\Omega} |\nabla w |^{p} dx \right)^{\frac{1}{p}},
\end{equation}
%\begin{align*}
%\displaystyle\int_{\Omega} |\nabla w |^{p} dx& \leq  \frac{c_{1}}{\alpha}(1+n)^{\theta+2}\int_{\Omega}  |\nabla w | dx \\
%&\leq c(n,\alpha) |\Omega |^{\frac{1}{p'} } \left(  \int_{\Omega} |\nabla w |^{p} dx  \right)^{\frac{1}{p}}, \\
%\end{align*}
we obtain
$$ \int_{\Omega} |\nabla w |^{p} dx  \leq c^{p^{\prime}}(n,\alpha)  |\Omega|,  $$
using the Poincar\'{e} inequality on the left hand side
\begin{equation}\label{l3}  ||w||_{L^{p}(\Omega)}\leq  c^{\frac{p'}{p}}  |\Omega |^{\frac{1}{p}}=c(n,\alpha,|\Omega |), \end{equation}
where $c(n,\alpha,|\Omega |)$ is a positive constant  independent form $v$, thus, we have that the ball $S$ of radius $c(n,\alpha,|\Omega |)$ is invariant for $G$.\\
Now, we are going to prove that the map $G$ is continuous in $S$. Consider a sequence $(v_{k})$ that converges to $v$ in $L^{p}(\Omega)$. We recall that $w_{k}=f_{n}h_{n}(v_{k})$  are bounded, we have that $w_{k}\in L^{\infty}(\Omega) $ and there exists a positive constant $d$, independent of $v_{k}$ and $w_{k}$, such that $||w_{k}||_{L^{\infty}(\Omega)}\leq d.$
Then by dominated convergence theorem
$$ ||f_{n}h_{n}(v_{k})- f_{n}h_{n}(v)||_{L^{p}(\Omega)}  \longrightarrow 0.$$
Hence, by the uniqueness of the weak solution,
we can say that $w_{k}=G(v_{k})$ converges to $w=G(v)$ in $L^{p}(\Omega)$. Thus $G $ is continuous over $L^{p}(\Omega)$.\\
What finally needs to be checked is that
$G(S)$ is relatively compact in $L^{p}(\Omega)$. Let $v_{k}$ be a bounded sequence, and let $w_{k}=G(v_{k}).$
Reasoning as to obtain \eqref{l3},  we have
$$\displaystyle\int_{\Omega} |\nabla w_{k}|^{p}dx =\displaystyle\int_{\Omega} |\nabla G(v_{k})|^{p}dx\leq c(n,\alpha,\gamma),$$
where $c$ is clearly independent from $v_{k}$ , so that, $G(L^{p}(\Omega))$ is relatively compact in  $L^{p}(\Omega)$.
Now, applying the Schauder's fixed point theorem that $G$ has a fixed point  $u_{n}\in S $ that is solution to \eqref{A1} in $W_{0}^{1,p}(\Omega)\cap L^{\infty}(\Omega).$\\
To show \eqref{l01}, we will consider the cases $m>1$ and $m=1.$\\
	Case $m>1,$ choosing $\varphi=|u_{n}|^{s(m-1)+\gamma} $ in \eqref{l02}, we have
	$$\displaystyle\int_{\Omega} |u_{n}|^{sm+\gamma}dx \leq \displaystyle\int_{\Omega} |f||u_{n}|^{s(m-1)}dx,$$
	therefore
	$$\displaystyle\int_{\Omega}|u_{n}|^{sm+\gamma} \leq c\left(\displaystyle\int_{\Omega}|f|^{m}\right)^{\frac{1}{m}}\left(\displaystyle\int_{\Omega}|u_{n}|^{sm+\gamma}\right)^{1-\frac{1}{m}},$$
	wich implies \eqref{l01}.\\
\par Case $m=1$. Choosing $\varphi=u_{n}^{\gamma}$, then
$$\displaystyle\int_{\Omega} |u_{n}|^{s-1} u_{n} u_{n}^{\gamma}dx \leq \displaystyle\int_{\Omega} \frac{f}{u_{n}^{\gamma}} u_{n}^{\gamma}dx\leq fdx,$$
which the estimate \eqref{l01}, as desired.
\begin{lem}
	Let $u_{n}$ be a solution to problem \eqref{A1} and $f\in L^{m}(\Omega)$ with $m\geq 1.$ Then
	\[
	\displaystyle\int_{\left\{k<u_{n}\right\}} u_{n}^{sm} \leq \frac{1}{k^{\theta}} \displaystyle\int_{\left\{k<u_{n}\right\}} f^{m}
	\]
	\[
	\operatorname{and} \displaystyle\lim _{|E| \rightarrow 0} \displaystyle\int_{E} u_{n}^{sm}=0 \text { uniformly with respect to } n.
	\]
\end{lem}
\textbf{Proof.} Let $k>0$ and $u_{n}^{s(m-1)}\psi_{i}$ be a sequence of increasing, positive, uniformly bounded $C^{\infty}(\Omega)$ functions, such that
\[
\psi_{i}(s) \rightarrow\left\{\begin{array}{ll}
1, & s \geq k \\
0, & 0 \leq s<k
\end{array}\right.
\]
The limit on $i$ gives
\[
\displaystyle\int_{\left\{k<u_{n}\right\}} u_{n}^{sm} \leq \displaystyle\int_{\left\{k<u_{n}\right\}} \frac{f u_{n}^{s(m-1)}}{\left(u_{n}+\frac{1}{n}\right)^{\theta}}
\]
Therefore we have
$$
\displaystyle\int_{\left\{k<u_{n}\right\}} u_{n}^{sm} \leq \frac{1}{\left(k+\frac{1}{n}\right)^{\theta}} \left( \displaystyle\int_{\left\{k<u_{n}\right\}} f^{m}\right)^{\frac{1}{m}}\left( \displaystyle\int_{\left\{k<u_{n}\right\}} u_{n}^{sm}\right) ^{1-\frac{1}{m}}$$
by Hölder inequality
\[
\displaystyle\int_{\left\{k<u_{n}\right\}} u_{n}^{sm} \leq \frac{1}{\left(k+\frac{1}{n}\right)^{\theta}} \int_{\left\{k<u_{n}\right\}} f^{m}
\]
This implies that
\[
\displaystyle\int_{E} u_{n}^{sm} \leq k^{r}|E|+\int_{E \cap\left\{u_{n}>k\right\}} u_{n}^{sm} \leq k^{r}|E|+\frac{1}{k^{\theta}} \displaystyle\int_{\left\{u_{n}>k\right\}} f^{m}
\]
since $f \in L^{m}(\Omega)$ for any given $\varepsilon>0,$ there exists $k_{\varepsilon}$ such that $\displaystyle\int_{\left\{\left|u_{n}\right|>k_{\varepsilon}\right\}}|f|^{m} \leq$
$\varepsilon .$ Therefore
\[
\displaystyle\int_{E} u_{n}^{sm} \leq k_{\varepsilon}^{sm}|E|+\frac{\varepsilon}{k^{\theta}}
\]
and the statement of this lemma is thus proved.
\begin{lem} \label{lem2}
	Let $u$ be a measurable function in $M^{r}(\Omega)$, $s>0,$ and suppose that there exists a positive constant $\rho>0$
	such that
	$$\displaystyle\int_{\Omega} |\nabla T_{k}(u)|^{p}dx\leq Ck^{\rho} \,\,\,\,\,\,\, \forall k>0.$$
	Then $|\nabla u|\in M^{\frac{p r}{\rho+r}}(\Omega).$
\end{lem}
\textbf{Proof}. Let $\lambda$ be fixed positive real number. For every $k>0,$ we have
$$meas\{|\nabla u|>\lambda\}=meas\{|\nabla u|> \lambda,|u|\leq k\}+meas\{|\nabla u|> \lambda,|u|> k\}$$
$$ \leq meas\{|\nabla u|> \lambda,|u|\leq k\}+meas\{|u|>k\}$$
and
$$meas\{|\nabla u|> \lambda,|u|\leq k\}\leq \frac{1}{\lambda^{p}}\int_{\Omega}|\nabla T_{k}(u)|^{p}dx\leq C\frac{k^{\rho}}{\lambda^{p}}.$$
Since $u\in M^{r}(\Omega)$, it follows that
$$meas\{|\nabla u|>\lambda\}=C\frac{k^{\rho}}{\lambda^{p}}+\frac{C}{k^{r}},$$
and this latter inequality holds for every $k>0$. Minimizing with respect to $k,$ we easily obtain
$$meas\{|\nabla u|>\lambda\}=\frac{C}{\lambda^{\frac{pr}{\rho+r}}}.$$
Thus, $|\nabla u|\in M^{\frac{p r}{\rho+r}}(\Omega).$
\begin{lem} \label{lem3}
	Let $u_{n}$ be a sequence of measurable functions such that $T_{k}(u_{n})$ is bounded in $W_{0}^{1,p}(\Omega)$
	for every $k>0$. Then there exists a measurable function $u$ such that $T_{k}(u)\in W_{0}^{1,p}(\Omega)$ and, moreover,
	$$T_{k}(u_{n})\longrightarrow  T_{k}(u)\,\,\,\,\mbox{weakly in}\,\,W_{0}^{1,p}(\Omega)\,\,\,\, \mbox{and}\,\,\,\,u_{n}\longrightarrow u\,\,\,a.e.\mbox{in}\,\,\,\Omega.$$
\end{lem}
\textbf{Proof}. Let us prove that $u_{n}\longrightarrow u$ locally in measure. To begin with, we observe that, for
$t,\varepsilon >0$, we have
$$\{|u_{n}-u_{m}|>t\}\subset \{|u_{n}|>k\} \cup \{|u_{m}|>k\}\cup \{|T_{k}(u_{n})-T_{k}(u_{m})|>t\}. $$
Therefore,
$$meas\{|u_{n}-u_{m}|>t\}\leq meas\{|u_{n}|>k\}+meas\{|u_{m}|>k\}+meas\{|T_{k}(u_{n})-T_{k}(u_{m})|>t\}. $$
Choosing $k$ large enough, we obtain
$$meas\{|u_{n}|>k\}<\varepsilon \,\,\,\, \mbox{and}\,\,\,meas\{|u_{m}|>k\}<\varepsilon.$$
We can assume that $\{T_{k}(u_{n})\}$ is a Cauchy sequence in $L^{q}(\Omega)$ for every\\ $q<p^{\ast}=\frac{Np}{N-p}.$ Then $$meas\{|T_{k}(u_{n})-T_{k}(u_{m})|>t\}\leq t^{-q}\int_{\Omega}|T_{k}(u_{n})-T_{k}(u_{m})|^{q}dx\leq \varepsilon \,\,\,\,\forall n,m\geq n_{0}(k,t).$$
This proves that $\{u_{n}\}$ is a Cauchy sequence in measure in $\Omega$. Therefore, there exists a measurable
function $u$ such that $u_{n}\longrightarrow u$ in measure. Hence that $u_{n}\longrightarrow u$ a.e.in $\Omega$, and so
$$T_{k}(u_{n})\longrightarrow T_{k}(u)\,\,\,\,\,\,\,\mbox{weakly in }\,\,\,W_{0}^{1,p}(\Omega).$$
%%%%%%%%%%%%%%%%%%%%%%%%%%%%%%%%%%%%%%%%%%%%%%%%%%%%%%%%%%%%%%%%%%%%%%%%%%%%%%%%%%%%%%%%%%%%%%%%%%%%%%%%%%%%%%%%%%%%%%%%%%%%%%%%%%%%%%%%%%%%%%%%%%%%%%%%%%%%%%%%%%%%%%%%%%%%%%%%%%%%%%%%%%%%%%%%%%%%%%%%%%%%%%%%%%%%%%%%%%%%%%%%%%%%%%%%%%%%%%%%%%%%%%%%%%%%%%%%%%%%lemmmeee%%%%%%%%%%%%%%%%%%%%%%%%%%%%%%%%%%%%%%%%%%%%%%%%%%%%%%%%%%%%%%%%%%%%%%%%%%%%%%
%	\begin{lem} \label{lem4}
%	Let $u_{n}$ is a solution to \eqref{A1}, such that $u_{n}$ is bounded in $W_{0}^{1,p}(\Omega)$ and satisfies %\eqref{l01} then
%	$$\int_{\Omega} h_{n}(u_{n})f_{n}\varphi\leq C$$
%and
% \begin{equation}\label{c}
% \lim_{n\longrightarrow+\infty} \int_{\Omega} h_{n}(u_{n})
% f_{n}\varphi dx =\int h(u)f\varphi dx,
%\end{equation}
%for all
%$\varphi \in W_{0}^{1,p}(\Omega)\cap L^{\infty} (\Omega).$ In both cases c is a positive constant independ of $n$ and $p$.
%	%\end{lem}
%%%%%%%%%%%%%%%%%%%%%%%%%%%%%%%%%%%%%%%%%%%%%%%%%%%%%%%%%%%%%%%%%%%%%%%%%%%%%%%%%%%%%%%%%%%%%%%%%%%%%%%%%%%%%%%%%%%%%%
\section {Proof of the results}\label{sec4}
This section is devoted to proving theorems cited above. We start with

%%%%%%%%%%%%%%%%%%%%%%%%%%%%%%%%%%%%%%%%%%%%%%%%%%%%%%%%%%%%%%%%%%%%%%%%%%%%%%%%%%%%%%%%%%%%%%%%%%%%%%%%%%%%%%%%%%%%%%%%%%%%%%%%%%%%%%%%%%%%%ùthoerms%%%%%%%%%%%%%%%%%%%%%%%%%%%%%%%%%%%%%%%%%%%%%%%%%%%%ùùùù

\textbf{Proof of Theorem \ref{thm1}.}
We separate our proof in three parts, according to the values of $s$ \\
\textbf{Part I.} Let $s\geq \frac{1+\theta-\gamma}{m-1}$ choosing  $\varphi=(1+u_{n})^{1+\theta}-1$ in \eqref{l02}, then, the second term is non negative. By \eqref{S2}, we have
$$\alpha\int_{\Omega}|\nabla u_{n}|^{p}dx \leq \int_{\Omega} |f| |u_{n}|^{1+\theta-\gamma} dx.$$
H\"{o}lder's inequality applied to the right-hand side yields
$$\int_{\Omega} |f| |u_{n}|^{1+\theta-\gamma} dx\leq c \left[ \displaystyle\int_{\Omega} u_{n}^{\frac{m(1+\theta-\gamma)}{m-1}}dx \right]^{1-\frac{1}{m}}. $$
So, we have
\begin{equation}\label{38}
\int_{\Omega} |\nabla u_{n}|^{p}dx \leq c\left[\displaystyle\int_{\Omega} u_{n}^{\frac{m(1+\theta-\gamma)}{m-1}}dx \right]^{1-\frac{1}{m}}.
\end{equation}
Since $$\frac{m(1+\theta-\gamma)}{m-1}\leq ms,$$
we have $s\geq \frac{1+\theta-\gamma}{m-1}$. Lemma \ref{lem1} implies that the right-hand side of \eqref{38} is uniformly bounded, so we have
\begin{equation}\label{d2}
\displaystyle\int_{\Omega} |\nabla u_{n}|^{p}dx \leq c.
\end{equation}
In order to prove that the limit function $u$ is a solution of \eqref{S1} in the sense of Definition$\ref{d1}$, we need to show that we can pass to the limit in the weak formulation of the approximating problems \eqref{A1}.\\
Now we focus on the left hand side of \eqref{l02}, by \eqref{d2} we conclude that there exist a subsequence, still indexed
by $n$, and a measurable function $u$ in $W_{0}^{1,p}(\Omega),$ such that
$u_{n}\rightharpoonup u$ weakly in $W_{0}^{1,p}(\Omega)$ and $u_{n}\longrightarrow u \,\,a.e \,\,\,in \,\, \Omega $.\\ Fatou's lemma implies $u \in L^{sm+\gamma} (\Omega).$ We see that
(see [\cite{04}, Lemma 5])
\begin{equation}\label{39}
\nabla u_{n} \longrightarrow \nabla u \,\,\,a.e\,\, in \,\,\,\Omega.
\end{equation}
Next, we pass to the limit in \eqref{l02}. By \eqref{39}, we can easily obtain
$$ |\nabla u_{n}|^{p-2}.|\nabla u_{n}| \rightharpoonup |\nabla u|^{p-2}.|\nabla u| \,\,\,weakly \,\,in \,\,  L^{p'}(\Omega).$$
Moreover,
$$ a(x,T_{n}(u_{n}))  \nabla \varphi \longrightarrow a(x,u) \nabla \varphi \,\,\,in\,\, L^{p}(\Omega).$$
Consequently, we have
$$ \int_{\Omega} a(x,T_{n}(u_{n})) |\nabla u_{n}|^{p-2}.|\nabla u_{n}|\nabla \varphi dx \longrightarrow \int_{\Omega} a(x,u) |\nabla u|^{p-2}.|\nabla u|\nabla \varphi dx. $$ Therefore, we can pass to the limit in the first term of the left-hand side of \eqref{l02}. We will show that
\begin{equation}\label{g1}
|u_{n}|^{s-1}u_{n}\rightarrow |u|^{s-1}u\,\,\,\,\mbox{in}\,\,\,L^{1}(\Omega).
\end{equation}
We take $ S_{\eta,k}(u_{n})$ as a test function in the weak formulation \eqref{A1}, we deduce
\begin{align*}
&\displaystyle\int_{\Omega} a(x,T_{n}(u_{n})) |\nabla u_{n}|^{p} S^{\prime}_{\eta,k}(u_{n})dx+\int_{\Omega}|u_{n}|^{s-1}u_{n}S_{\eta,k}(u_{n})dx \\
&\leq \sup_{s\in [k,\infty)}[h(s)]\int_{\Omega}f_{n}S_{\eta,k}(u_{n}),\\
\end{align*}
which, observing that the first term on the left hand side is non negative and taking the limit with respect to
$\eta \rightarrow 0,$ implies
$$\displaystyle\int_{\{u_{n}\geq k\}} |u_{n}|^{s-1}u_{n}dx\leq \sup_{s\in [k,\infty)}[h(s)]\int_{\{u_{n}\geq k\}} f_{n}dx,$$
which, since $f_{n}$ converges to $f$ in $L^{m}(\Omega)$, easily implies that $|u_{n}|^{s-1}u_{n} $ is equi-integrable and so it converges to $|u|^{s-1}u$ in $L^{1}(\Omega),$ this concludes \eqref{g1}.\\
The next step we want to pass to the limit in the right hand side of \eqref{l02}.
Let us take $0\leq \varphi\in W_{0}^{1,p}(\Omega)\cap L^{\infty}(\Omega) $ as test function in the weak formulation of \eqref{A1}, by using the young inequality and the hypotheses in \eqref{S2} and \eqref{S3}, we have
\begin{align*}
&\displaystyle\int_{\Omega} h_{n}(u_{n}) f_{n}  \varphi =
\int_{\Omega} a(x,T_{n}(u_{n})) |\nabla u_{n}|^{p-2}.|\nabla u_{n}|\nabla \varphi dx+\int_{\Omega}u_{n}^{s-1}u_{n}\varphi dx\\
& \leq C||\varphi||_{L^{\infty}(\Omega)} +\beta \displaystyle\int_{\Omega}|\nabla u_{n}|^{p-1} \nabla \varphi dx \\
& \leq C||\varphi||_{L^{\infty}(\Omega)}
+\beta\frac{p-1}{p}\displaystyle\int_{\Omega}|\nabla u_{n}|^{p} dx+  \beta\frac{1}{p}\displaystyle\int_{\Omega}|\nabla \varphi |^{p} dx \\
& \leq C||\varphi||_{L^{\infty}(\Omega)}+C\left[ \displaystyle\int_{\Omega}|\nabla \varphi |^{p}dx+\displaystyle \displaystyle\int_{\Omega}|\nabla u_{n}|^{p}dx \right],\\
\end{align*}
then
\begin{equation}\label{41}
\int_{\Omega} h_{n}(u_{n}) f_{n}  \varphi \leq C||\varphi||_{L^{\infty}(\Omega)}+C [||\varphi||_{W_{0}^{1,p}(\Omega)}+||u_{n}||_{W_{0}^{1,p}(\Omega)}].
\end{equation}
From now on, we assume that $h(s)$ is unbounded as $s$ tends to $0.$
%%	From now we consider a non negative $ \varphi \in W_{0}^{1,p}(\Omega) \cap L^{\infty}(\Omega).$
An application of the Fatou Lemma in \eqref{41}	
with respect to $ n $ gives
\begin{equation}\label{42}
\int_{\Omega} h(u)f \varphi \leq c,
\end{equation}
where $c$ does not depend on $n$.\\
Hence $fh(u)\varphi \in L^{1}(\Omega)$ for any nonnegative $\varphi\in W^{1,p}_{0}(\Omega)$. As a consequence, if $h(s)$ is unbounded as $s$
tends to $ 0$, we deduce that
\begin{equation}\label{u0}
\{u=0\}\subset \{f=0\}
\end{equation}
up to a set of zero Lebesgue measure.\\
Now, for $\delta>0,$ we split the right hand side of \eqref{l02} as
\begin{equation}\label{u1}
\displaystyle\int_{\Omega} h_{n}(u_{n}) f_{n}  \varphi dx =\displaystyle\int _{\{u_{n}\leq \delta\}} h_{n}(u_{n}) f_{n}  \varphi dx+
\displaystyle\int_{\{u_{n}> \delta\}} h_{n}(u_{n}) f_{n} \varphi dx,
\end{equation}
and we pass to limit as $ n\rightarrow +\infty$ and then $\delta \rightarrow 0 $, we remark that we need to choose
$\delta \neq \{\eta ;  |u=\eta|> 0\},$
which is at most a countable set, for the second term \eqref{u1} we have
\begin{equation}\label{314}
0\leq	h_{n}(u_{n}) f_{n}  \varphi \chi_{\{u_{n}>\delta\}} \leq \sup_{s\in ]\delta,\infty)} [h(s)]f \varphi \in L^{1}(\Omega),
\end{equation}
which precis to apply the Lebesgue Theorem  with respect $n$. Hence on has
$$\lim_{n \rightarrow +\infty}\int _{\{u_{n} > \delta\}} h_{n}(u_{n}) f_{n}  \varphi dx= \int _{\{u > \delta\}} h(u) f \varphi dx.$$
Moreover it follows by \eqref{42}that
\begin{equation}\label{14}
\displaystyle\lim_{\delta \rightarrow 0^{+}}\displaystyle\lim_{n \rightarrow +\infty}\displaystyle\int _{\{u_{n} > \delta\}} h_{n}(u_{n}) f_{n}  \varphi dx
=\displaystyle\int _{\{u > 0\}} h(u) f \varphi dx.
\end{equation}
Now in order to get rid of the first term of the right hand side of \eqref{u1}, we take $ V_{\delta}(u_{n})\varphi$ is a test function in the weak formulation of \eqref{A1},where\\ $V_{\delta}(u_{n}):=V_{\delta,\delta}(u_{n})$ is defined in \eqref{N1} and by Lemma 1.1 contained in \cite{03}, we have $V_{\delta}(u_{n})$  belongs to $W_{0}^{1,p}(\Omega),$ then (recall $V^{'}_{\delta}(u_{n})\leq 0\,\, for\,\, s\geq 0)$
\begin{align*}
&\displaystyle\int _{\{u_{n} \leq \delta\}} h_{n}(u_{n}) f_{n}  \varphi dx\leq \displaystyle\int _{\Omega} h_{n}(u_{n}) f_{n} V_{\delta}(u_{n}) \varphi dx\\
&=\displaystyle\int_{\Omega} a(x,T_{n}(u_{n}))|\nabla u_{n}|^{p-2}\nabla u_{n}\nabla \varphi V_{\delta}(u_{n})dx\\
&-\frac{1}{\delta}\displaystyle\int_{\{\delta<u_{n}<2\delta\}} a(x,T_{n}(u_{n}))|\nabla u_{n}|^{p-2}\nabla u_{n} \varphi \nabla u_{n}dx
+\displaystyle\int_{\Omega}|u_{n}|^{s-1}u_{n}V_{\delta}(u_{n})\varphi dx,\\
\end{align*}
by using \eqref{S2} and \eqref{S3}, we have \\
\begin{align*}
\displaystyle\int _{\{u_{n} \leq \delta\}} h_{n}(u_{n}) f_{n}  \varphi dx&\leq \beta \int_{\Omega}|\nabla u_{n}|^{p-2}\nabla u_{n}\nabla \varphi V_{\delta}(u_{n})dx\\
&+\int_{\Omega}|u_{n}|^{s-1}u_{n}V_{\delta}(u_{n}) \varphi dx,\\
\end{align*}
using that $V_{\delta}$ is bounded we deduce that $|\nabla u_{n}|^{p-2}\nabla u_{n} V_{\delta}(u_{n})$ converges to
$|\nabla u|^{p-2}\nabla u V_{\delta}(u)$ weakly in $L^{p^{\prime}}(\Omega)^{N}$ as $n$ tends to infinity. This implies that
\begin{equation}\label{12}
\lim_{n\rightarrow +\infty}\int _{\{u_{n} \leq \delta\}} h_{n}(u_{n}) f_{n}  \varphi dx\leq \beta \int_{\Omega}|\nabla u|^{p-2}\nabla u\nabla \varphi V_{\delta}(u)dx+\int_{\Omega}|u|^{s-1}uV_{\delta}(u) \varphi dx.
\end{equation}
Since $V_{\delta}(u)$ converges to $\chi_{\{u=0\}}$ a.e in $\Omega$ as $\delta$ tends to $0$ and since
$u \in W_{0}^{1,p}(\Omega),$ then $|\nabla u|^{p-2}\nabla u\nabla \varphi V_{\delta}(u)$ converges to $0$ a.e.
in 	$\Omega$ as $\delta$ tends to $0.$ Applying the Lebesgue Theorem on the right hand side of \eqref{12} we obtain that
\begin{equation}\label{13}
\lim_{\delta\rightarrow 0^{+}} \lim_{n\rightarrow +\infty}\int _{\{u_{n} \leq \delta\}} h_{n}(u_{n}) f_{n}  \varphi dx\leq
\beta \int_{\{u=0\}}|\nabla u|^{p-2}\nabla u\nabla \varphi dx+\int_{\{u=0\}}|u|^{s-1}u\varphi dx
=0,
\end{equation}	
by \eqref{14} and \eqref{13}, we deduce that
\begin{equation}\label{44}
\displaystyle\lim_{n \rightarrow +\infty}\int_{\Omega} h_{n}(u_{n}) f_{n}  \varphi dx
=\displaystyle\int_{\Omega} h(u) f \varphi dx\,\,\,\,\forall 0\leq \varphi \in W_{0}^{1,p}(\Omega)\cap L^{\infty}(\Omega).
\end{equation}
Moreover, decomposing any $\varphi=\varphi^{+}-\varphi^{-}$, and using that \eqref{44} is linear in $\varphi$, we deduce that \eqref{44} holds for every
$\varphi \in W_{0}^{1,p}(\Omega)\cap L^{\infty}(\Omega)$.\\
We treated $h(s)$ unbounded as $s$ tends to $0$, as regards bounded function $h$ the proof is easier and only difference
deals with the passage to the limit in the left hand side of \eqref{44}. We can avoid introducing $\delta$ and we can  substitute \eqref{314} with
$$0\leq f_{n}h_{n}(u_{n})\varphi \leq f||h||_{L^{\infty}(\Omega)}\varphi.$$
Using the same argument above we have that
\begin{equation}\label{cn}
\displaystyle\lim_{n\rightarrow +\infty} \int_{\Omega}f_{n}h_{n}(u_{n})\varphi dx=\displaystyle\int_{\Omega}fh(u)\varphi dx,
\end{equation}
whence one deduces \eqref{S5}. This concludes the proof of part I.\\
\textbf{Part II.} Let $$\frac{1+\theta -\gamma}{pm-1}<s< \frac{1+\theta-\gamma}{m-1}.$$
Choosing $\varphi =(1+u_{n})^{s(m-1)+\gamma}-1 $ in \eqref{l02}, we  see that the second term is non negative. Using assumption \eqref{S2}, we have
$$\displaystyle\int_{\Omega} \frac{|\nabla u_{n}|^{p}}{(1+u_{n})^{1+\theta-s(m-1)-\gamma}}dx \leq c \displaystyle\int_{\Omega} |f| |u_{n}|^{s(m-1)} dx.$$
Now, using H\"{o}lder's inequality in the right-hand side of the previous inequality, we obtain
$$\displaystyle\int_{\Omega} |f| |u_{n}|^{s(m-1)} dx\leq c \left[\displaystyle\int_{\Omega} u_{n}^{ms+\gamma} dx \right]^{1-\frac{1}{m}} \leq c .$$
Therefore,
\begin{equation}\label{310}
\displaystyle\int_{\Omega} \frac{|\nabla u_{n}|^{p}}{(1+u_{n})^{1+\theta-s(m-1)-\gamma}}dx \leq c.
\end{equation}
Let $1\leq \sigma < p.$ Let us write
$$\displaystyle\int_{\Omega} |\nabla u_{n}|^{\sigma} dx=\int_{\Omega} \frac{|\nabla u_{n}|^{\sigma}(1+u_{n})^{\frac{\sigma}{p}(1+\theta-s(m-1)-\gamma)}}{(1+u_{n})^{\frac{\sigma}{p}(1+\theta-s(m-1)-\gamma)}}dx $$
H\"{o}lder's inequality implies
$$
\displaystyle\int_{\Omega} \frac{|\nabla u_{n}|^{\sigma}(1+u_{n})^{\frac{\sigma}{p}(1+\theta-s(m-1)-\gamma)}}{(1+u_{n})^{\frac{\sigma}{p}(1+\theta-s(m-1)-\gamma)}}dx \leq c \left[\int_{\Omega}(1+u_{n})^{\frac{\sigma}{p-\sigma}[1+\theta-s(m-1)-\gamma]}dx \right]^{1-\frac{\sigma}{p}}.$$
Therefore, we have
\begin{equation}\label{311}
\displaystyle\int_{\Omega} |\nabla u_{n}|^{\sigma}dx \leq c \left[\displaystyle\int_{\Omega}(1+u_{n})^{\frac{\sigma}{p-\sigma}[1+\theta-s(m-1)-\gamma]}dx \right]^{1-\frac{\sigma}{p}}.
\end{equation}
If $$ \frac{\sigma}{p-\sigma}[1+\theta -s(m-1)-\gamma]\leq ms \,\,i.e\,\,\sigma \leq \frac{pms}{s+1+\theta-\gamma},$$
then Lemma \ref{lem1} implies that the right- hand side of \eqref{311} is uniformly bounded. The inequality
$$\frac{1+\theta-\gamma}{pm-1}<s \,\,\mbox{ implies } \,\,\frac{pms}{s+1+\theta-\gamma}>1. $$
Consequently,
$$\displaystyle\int_{\Omega} |\nabla u_{n}|^{\sigma}dx \leq c,\,\, \sigma=\frac{pms}{1+\theta+s-\gamma}.$$
Up to a subsequence, there exists a function $u\in W_{0}^{1,\sigma}(\Omega)$ such that
$$ u_{n}\rightharpoonup u\,\, \mbox{weakly} \,\,in \, W_{0}^{1,\sigma}(\Omega)\,\,and\,\, u_{n}\longrightarrow u\,\, a.e\,\,in\,\,\,\Omega. $$
By Lemma 5 (see\cite{04}), we have $ \nabla u_{n}\longrightarrow \nabla u \,\,\,a.e\,\,\,in \,\,\Omega.$ Fatou's Lemma implies $u_{n}^{sm+\gamma} \in L^{1}(\Omega)$ we will now pass to the limit in \eqref{l02}. We can easily obtain
$$|\nabla u_{n}|^{p-2} |\nabla u_{n}| \longrightarrow  |\nabla u|^{p-2} \nabla u\,\,\,weakly \,\,in\,\,\, L^{\frac{\sigma}{p-1}}(\Omega),$$
and
$$ a(x,T_{n}(u_{n})) \varphi  \longrightarrow  a(x,u) \varphi,\,\,in\,\,L^{(\frac{\sigma}{p-1})^{\prime}}(\Omega).$$
Therefore,we have
$$ \displaystyle\int_{\Omega}a(x,T_{n}(u_{n}))|\nabla u_{n}|^{p-2} \nabla u_{n} \nabla \varphi dx \longrightarrow \int_{\Omega} a(x,u)|\nabla u|^{p-2} \nabla u \nabla \varphi dx.$$
The remaining two parts in \eqref{l02} are the same as part I.\\
PART III. Let $0< s \leq \frac{1+\theta -\gamma}{pm-1}.$ It follows from \eqref{310} that
$$\displaystyle\int_{\Omega \cap\{|u_{n}|<k \}}\frac{|\nabla u_{n}|^{p}}{(1+u_{n})^{1+\theta-s(m-1)-\gamma}}dx \leq c,$$ and consequently
\begin{equation}\label{312}
\displaystyle\int_{\Omega}|\nabla T_{k}(u_{n})|^{p}dx =\displaystyle\int_{\Omega \cap\{|u_{n}|<k \}}|\nabla T_{k}(u_{n})|^{p}dx\leq c(1+k)^{1+\theta-s(m-1)-\gamma}.
\end{equation}
Lemma$\ref{lem3}$ implies the existence of a measurable function $u$ such that \\$T_{k}(u)\in W_{0}^{1,p}(\Omega)$
for any $k>0$, besides, passing if necessary to subsequence, we have
$$T_{k}(u_{n})\rightharpoonup T_{k}(u)\,\,\,\mbox{weakly in } W_{0}^{1,p}(\Omega)\,\,\,\mbox{and}\,\,\,a.e.\mbox{in}\,\, \Omega.$$
Fatou's Lemma implies that $|u|^{s}\in L^{1}(\Omega)$. We can pass to the limit in \eqref{312}, to get
$$\int_{\Omega}|\nabla T_{k}(u)|^{p}dx \leq c(1+k)^{1+\theta-s(m-1)-\gamma}.$$
Since
$$s\leq \frac{1+\theta -\gamma}{pm-1}\leq \frac{1+\theta -\gamma}{m-1},	$$
we have $1+\theta-s(m-1)-\gamma >0 $.\\
As a result of the Lemma \ref{lem2}, we obtain $| \nabla u| \in M^{\frac{pms}{1+\theta+s-\gamma}}(\Omega).$
We will show that $u$ is an entropy solution of \eqref{S1}. Indeed, let us choose
$$T_{k}(u_{n}-\varphi),\,\,\,\,\,\varphi\in	W_{0}^{1,p}(\Omega)\cap L^{\infty}(\Omega),$$
as a test function in \eqref{l02}, then we have $$
\displaystyle\int_{\Omega} a(x,T_{n}(u_{n}))|\nabla u_{n}|^{p-2}\nabla u_{n}\nabla T_{k}(u_{n}-\varphi) dx + \int_{\Omega}|u_{n}|^{s-1} u_{n} T_{k}(u_{n}-\varphi)$$
\begin{equation}\label{k1}
=  \displaystyle\int_{\Omega}f_{n}h_{n}(u_{n})T_{k}(u_{n}-\varphi).
\end{equation}
Let us pass to the limit in \eqref{k1}. For the second term on the left-hand side and for the right-hand side, we can use \eqref{cn}to obtain the limit. For the first term on the left-hand side, we will firstly show that
$\nabla T_{k}(u_{n})\rightarrow \nabla T_{k}(u)$ a.e. in  $\Omega$. Let $ \varphi=T_{k}(u_{n})-T_{k}(u)$ in \eqref{l02}, then we obtain
%		$$T_{k}(u_{n}-\varphi),\,\,\,\,\,\varphi\in	W_{0}^{1,p}(\Omega)\cap L^{\infty}(\Omega)$$
%%	as a test function in \eqref{l02}, then we have $$
$$	\int_{\Omega} a(x,T_{n}(T_{k}(u_{n})))|\nabla u_{n}|^{p-2}\nabla u_{n}[\nabla	T_{k}(u_{n})-\nabla T_{k}(u)]dx + \int_{\Omega}|u_{n}|^{s-1} u_{n} [T_{k}(u_{n})-T_{k}(u)]$$
$$=  \int_{\Omega}f_{n}h_{n}(u_{n})[T_{k}(u_{n})-T_{k}(u)].$$
As a consequence, we have
$$\int_{\Omega} a(x,T_{n}(T_{k}(u_{n})))
[|\nabla T_{k}(u_{n})|^{p-2}\nabla T_{k}(u_{n})-|\nabla T_{k}(u)|^{p-2}\nabla T_{k}(u)]
[\nabla	T_{k}(u_{n})-\nabla T_{k}(u)]dx $$
$$=  \int_{\Omega}f_{n}h_{n}(u_{n})[T_{k}(u_{n})-T_{k}(u)]dx- \int_{\Omega}|u_{n}|^{s-1} u_{n} [T_{k}(u_{n})-T_{k}(u)]dx$$	
\begin{equation}\label{k2}
-\int_{\Omega} a(x,T_{n}(T_{k}(u_{n})))|\nabla T_{k}(u)|^{p-2}\nabla T_{k}(u)]
[\nabla	T_{k}(u_{n})-\nabla T_{k}(u)]dx.
\end{equation}
We are going to show that the three terms of the right-hand side in \eqref{k2} all converge to zero. For the first term, we can use the \eqref{cn} to take the limit. As the result of the proof in part one, we obtain
$$ u_{n}^{s-1}u_{n} \rightarrow |u|^{s-1}u	\,\,\, \mbox{in} \,\,\, L^{1}(\Omega).$$
Therefore, we have
$$ \int_{\Omega}  u_{n}^{s-1}u_{n}[T_{k}(u_{n})-T_{k}(u)]dx \rightarrow	 0 \,\,\,\mbox{as}\,\,\,n\rightarrow \infty.$$
We can easily know the fact that $a(x,T_{n}(T_{k}(u_{n})))|\nabla T_{k}(u)|^{p-2}\nabla T_{k}(u)\in L^{p^{\prime}}(\Omega).$ Thus, for every measurable set $E \subset \Omega$, we can write
$$ \int_{E} |a(x,T_{n}(T_{k}(u_{n})))|\nabla T_{k}(u)|^{p-1}|^{p^{\prime}}dx\rightarrow	 0 \,\mbox{as}\, meas E\rightarrow 0.$$
Because
$$a(x,T_{n}(T_{k}(u_{n})))|\nabla T_{k}(u)|^{p-2}\nabla T_{k}(u)\rightarrow a(x,T_{k}(u))|\nabla T_{k}(u)|^{p-2}\nabla T_{k}(u)\, \mbox{a.e.in}\,\Omega,$$
we have
$$a(x,T_{n}(T_{k}(u_{n})))|\nabla T_{k}(u)|^{p-2}\nabla T_{k}(u)\rightarrow a(x,T_{k}(u))|\nabla T_{k}(u)|^{p-2}\nabla T_{k}(u)\,\mbox{in}\,\,L^{p^{\prime}}(\Omega),$$	
by Vitali's Theorem. By Lemma \ref{lem3}, we see that
$$\nabla T_{k}(u_{n})-\nabla T_{k}(u)\rightharpoonup 0\,\,\,\,\mbox{weakly in}\,\,L^{p}(\Omega).$$
Therefore,
$$\int_{\Omega} a(x,T_{n}(T_{k}(u_{n})))|\nabla T_{k}(u)|^{p-2}\nabla T_{k}(u)]
[\nabla	T_{k}(u_{n})-\nabla T_{k}(u)]dx \rightarrow 0 \,\,\,\mbox{as}\,\,\,n \rightarrow \infty.$$
From the above, we have
$$\int_{\Omega} a(x,T_{n}(T_{k}(u_{n})))|\nabla T_{k}(u)|^{p-2}\nabla T_{k}(u)]
[\nabla	T_{k}(u_{n})-\nabla T_{k}(u)]dx \rightarrow 0.$$
As a consequence, Lemma 5 in \cite{04} implies $ \nabla T_{k}(u_{n})\rightarrow\nabla T_{k}(u)$ in $L^{p}(\Omega)$. Therefore,
$$ \nabla T_{k}(u_{n})\rightarrow \nabla T_{k}(u)\,\,\,\mbox{a.e. in}\,\,\Omega.$$
Let $m=k+|\varphi|.$ The first term on the left-hand side in \eqref{k1} can be rewritten as
$$\int_{\Omega}a(x,T_{n}(u_{n}))|\nabla T_{m}(u)|^{p-2}\nabla T_{m}(u) \nabla T_{k}(u_{n}-\varphi)dx.$$
Since
$ \nabla T_{m}(u_{n})\rightarrow\nabla T_{m}(u)$ a.e.in $\Omega$, as a result of the Fatou's Lemma, we have
$$\liminf_{n\rightarrow \infty}\int_{\Omega}a(x,T_{n}(u_{n}))|\nabla T_{m}(u)|^{p-2}\nabla T_{m}(u) \nabla T_{k}(u_{n}-\varphi)dx$$
$$\geq \int_{\Omega}a(x,u)|\nabla T_{m}(u)|^{p-2}\nabla T_{m}(u) \nabla T_{k}(u_{n}-\varphi)dx$$
$$= \displaystyle\int_{\Omega}a(x,u)|\nabla u|^{p-2}\nabla u \nabla T_{k}(u_{n}-\varphi)dx.\,\,\,\,\,\,\,\,\,\,\,\,\,\,\,\,\,\,\,\,\,\,\,\,\,\,\,$$
So we see that $u$ is an entropy solution of \eqref{S1}.

\textbf{Proof of Theorem$\ref{t2}$.}
We separate our proof in two parts, according to the values of $s$ \\
\textbf{Part I.} Let $0<s\leq \frac{N(1-\gamma )+\gamma }{m(N-1)}.$ Then $s\leq \frac{N(1-\gamma )+\gamma }{m(N-1)},$ which implies $ms+\gamma \leq \sigma^{*}$ for $\sigma \geq 1$. As the result of \eqref{311} and of Sobolev's embedding theorem, we have
$$\displaystyle\int_{\Omega}  |u_{n}|^{\sigma^{\ast}} dx
\leq c \left[\displaystyle\int_{\Omega}(1+u_{n})^{\frac{\sigma}{p-\sigma}[1+\theta-s(m-1)-\gamma]}dx  \right]^{\frac{(p-\sigma)\sigma^{\ast}}{p\sigma}}.$$
If $$\frac{\sigma}{p-\sigma}[1+\theta-s(m-1)-\gamma]\leq \sigma^{\ast}\,\,ie\,\, \sigma \leq \frac{N[p+s(m-1)-1-\theta+\gamma]}{N+s(m-1)-1-\theta+\gamma},$$
then, by $m> 1$  and $p>p_{0}> 1+\frac{(N-1)[1+\theta-s(m-1)-\gamma]}{N},$
which implies
$$ \frac{N[p+s(m-1)-1-\theta+\gamma]}{N+s(m-1)-1-\theta+\gamma}>1.$$
\begin{equation}\label{eq1}
\int_{\Omega}  |u_{n}|^{\sigma^{\ast}} dx \leq c+c \left(\displaystyle\int_{\Omega}  |u_{n}|^{\sigma^{*}} dx\right)^{\frac{(p-\sigma)\sigma^{\ast}}{p\sigma}}.
\end{equation}
For \eqref{eq1}, by Young's inequality with $\epsilon $, we have
$$\int_{\Omega}  |u_{n}|^{\sigma^{\ast}} dx
\leq c. $$
Which, to gether \eqref{311} and
$\frac{\sigma}{p-\sigma}[1+\theta-s(m-1)-\gamma]\leq\sigma^{\ast},$
implies
$$ \int_{\Omega}  |\nabla u_{n}|^{\sigma}dx \leq c,\,\,\,\, \sigma \leq \frac{N[p+s(m-1)-1-\theta+\gamma]}{N+s(m-1)-1-\theta+\gamma}.$$
The remaining proof of this part is the same as part II in Theorem \ref{thm1}, we have can show that $u$ is a distributional solution to problem \eqref{S1}.\\
\textbf{Part II.} Let $s \geq \frac{N(1-\gamma )+\gamma }{m(N-1)}$. Since $p> p_{0}$, it follows that
$$\frac{N(1-\gamma )+\gamma }{m(N-1)} > \frac{1+\theta -\gamma}{pm-1},$$
thus, we can show that $u$ is a distributional solution to the problem \eqref{S1} by the same method as in Part II of Theorem \ref{thm1}.\\
\textbf{Proof of Theorem \ref{t3}.}
We separate our proof in two parts, according to the values of $s$ \\
\textbf{Part a.} Let $s> \frac{1+\theta -\gamma}{p-1}$. Choosing $\varphi =(1+u_{n})^{\gamma}-1$
in \eqref{l02}, then the second term is non negative. Using assumption \eqref{S2}, we can write
\begin{equation}\label{31}
\int_{\Omega}  \frac{|\nabla u_{n}|^{p}}{(1+u_{n})^{1+\theta -\gamma}} dx
\leq \int_{\Omega} |f| dx \leq c.
\end{equation}
Let $r<p,$ writing
$$
\displaystyle\int_{\Omega}  |\nabla u_{n}|^{r} dx
=\displaystyle\int_{\Omega}  \frac{|\nabla u_{n}|^{r}}{(1+u_{n})^{\frac{r(1+\theta -\gamma)}{p}}} (1+u_{n})^{\frac{r(1+\theta -\gamma)}{p}}dx. $$
$$\leq \left( \int_{\Omega}  \frac{|\nabla u_{n}|^{p}}{(1+u_{n})^{(1+\theta -\gamma)}} dx \right)^{\frac{r}{p}}
\left(\displaystyle\int_{\Omega}(1+u_{n})^{\frac{r(1+\theta -\gamma)}{p-r}}dx \right)^{1-\frac{r}{p}}
$$
$$\leq c \left(\displaystyle\int_{\Omega}(1+u_{n})^{\frac{r(1+\theta -\gamma)}{p-r}}dx \right)^{1-\frac{r}{p}}.$$
Thanks to Lemma \ref{lem1}, if
$$\frac{r}{p-r}(1+\theta-\gamma)\leq s\,\,_,\, ie\,\,\,\, r< \frac{ps}{1+\theta+s-\gamma}.$$
Then $$s>\frac{1+\theta-\gamma}{p-1}\,\,\,\,\mbox{ implies }\,_,\,\,\,\frac{ps}{1+\theta+s-\gamma}>1.$$
In that case, the right-hand sides is uniformly
bounded  and so we have
$$\displaystyle\int_{\Omega} |\nabla u_{n}|^{r}dx\leq c\,\,\,\,\,,r<\frac{ps}{1+s+\theta-\gamma}. $$
As a consequence, there exists a function $u\in W_{0}^{1,r}(\Omega)$ such that
$$ u_{n}\rightharpoonup u\,\,\, weakly \,\,in \,\, W_{0}^{1,r}(\Omega)\,\,and\,\, u_{n}\longrightarrow u\,\,\, a.e\,\,in\,\,\,\Omega. $$
Let
$$g_{n}=f_{n}h_{n}(u_{n})-T_{n}(|u_{n}|^{s-1}u_{n}).$$
Because $g_{n}$ is bounded in $L^{1}(\Omega),$ and $u_{n}$ is a solution of\\
\begin{equation*}
\left\{
\begin{array}{ll}
-\mbox{div} \big(a(x,T_{n}(u_{n}))|\nabla u_{n}|^{p-2}|\nabla u_{n}|\big)=g_{n}, \\
u_{n}\in W_{0}^{1,p}(\Omega),
\end{array}
\right.
\end{equation*}
it follows from Lemma 1 (see\cite{05}), that
\begin{equation}\label{32}
\nabla u_{n}\longrightarrow \nabla u \,\,\,a.e\,\,\,in \,\,\Omega.
\end{equation}
We are going to show that $u$ is a distributional solution to problem \eqref{S1} by passing  to the limit in \eqref{l02}. We suppose that $\varphi\in C_{0}^{\infty}(\Omega).$ Since \\$|\nabla u_{n}|^{p-2} |\nabla u_{n}|\in L^{\frac{r}{p-1}(\Omega)}$ and  \eqref{32} hold, we have
$$|\nabla u_{n}|^{p-2} |\nabla u_{n}| \longrightarrow  |\nabla u|^{p-2} \nabla u\,\,\,weakly \,\,in\,\,\, L^{\frac{r}{p-1}}(\Omega).$$
Vitali's Theorem implies that
$$ a(x,T_{n}(u_{n})) \nabla \varphi  \longrightarrow  a(x,u) \nabla \varphi,\,\,\,\,\,in\,\,L^{(\frac{r}{p-1})^{\prime}}(\Omega),$$
where $(\frac{r}{p-1})^{\prime}=\frac{p-1-r}{p-1}.$ Therefore, we can pass to the limit in the first term on the left-hand side of \eqref{l02}. For the second term on the left hand-side and the first term on the right-hand side in \eqref{l02} we can namely arguing exactly as part I in Theorem$\ref{thm1}$. Therefore, we conclude that $u$ is a distributional solution to problem \eqref{S1}.\\
\textbf{Part b.} Let $ 0<s\leq \frac{1+\theta-\gamma}{p-1}$. Let us choose $T_{k}(u_{n})$ as a test function in \eqref{l02}; then the second term is non negative. Using assumption \eqref{S2} we can write
\begin{equation}\label{34}
\int_{\Omega} |\nabla T_{k}(u_{n})|^{p}dx\leq c k^{1-\gamma}(1+k)^{\theta}\leq c(1+k)^{1-\gamma+\theta}.
\end{equation}
By Lemma \ref{lem3}, there exists a function $u$ such that
$T_{k}(u)\in W_{0}^{1,p}(\Omega)$. Moreover,
$$T_{k}(u_{n})\rightharpoonup  T_{k}(u)\,\,\,weakly \,\, in \,\, W_{0}^{1,p}(\Omega) \,\, \forall k>0 \,\,\,\, and\,\,\, u_{n }\longrightarrow u \,\, a.e\,\,\,in \,\, \Omega $$
Fatou's Lemma implies that $|u|^{s}\in L^{1}(\Omega). $ We can pass to the limit in \eqref{34},to get
$$\int_{\Omega}|\nabla T_{k}(u)|^{p}dx \leq c(1+k)^{1+\theta-\gamma}.$$
As a result of the Lemma$\ref{lem2}$, we obtain $| \nabla u| \in M^{\frac{ps}{1+\theta+s-\gamma}}(\Omega).$

By the same method as in part II of Theorem \ref{thm1}, we can show that $u$ is an entropy solution of \eqref{S1}.\\
\textbf{Proof of theorem $\ref{t4}$}.
We separate our proof into three parts, according to the values of $s$ \\
\textbf{Part 1.} Let $0<s\leq \frac{N(1-\gamma)+\gamma}{N-1}.$
It is obvious that  $s\leq \frac{N(1-\gamma)+\gamma}{N-1}$
implies $s+\gamma\leq r^{\ast}.$ Using \eqref{31}, we obtain
$$\int_{\Omega}  \frac{|\nabla u_{n}|^{p}}{(1+u_{n})^{1+\theta -\gamma}} dx
\leq c. $$
Let $ 1\leq r<p,$ let us write
$$
\displaystyle\int_{\Omega}  |\nabla u_{n}|^{r} dx
=\displaystyle\int_{\Omega}  \frac{|\nabla u_{n}|^{r}}{(1+u_{n})^{\frac{r(1+\theta -\gamma)}{p}}} (1+u_{n})^{\frac{r(1+\theta -\gamma)}{p}}dx.$$
By H\"{o}lder's inequality, we have
\begin{equation}\label{37}
\displaystyle\int_{\Omega}  |\nabla u_{n}|^{r} dx \leq
c 	\left(	\displaystyle\int_{\Omega}(1+u_{n})^{\frac{r(1+\theta -\gamma)}{p-r}}dx \right)^{1-\frac{r}{p}}.
\end{equation}
Sobolev's embedding Theorem implies
$$	\left(	\displaystyle\int_{\Omega}  | u_{n}|^{r^{\ast}} dx\right)^{\frac{1}{r^{\ast}}}\leq c \left(\int_{\Omega}  |\nabla u_{n}|^{r} dx\right)^{\frac{1}{r}}\,, r^{\ast}=\frac{Nr}{N-r}.$$
Therefore,
$$	\left(\displaystyle\int_{\Omega}  | u_{n}|^{r^{\ast}} dx \right)^{\frac{1}{r^{\ast}}} \leq c \left(\displaystyle\int_{\Omega}  (1+u_{n})^{\frac{r(1+\theta-\gamma)}{p-r}}dx \right)^{\frac{(p-r)r^{\ast}}{pr}}.$$
Suppose $\frac{r(1+\theta-\gamma)}{p-r} \leq r^{\ast}$ ie,$\,\,
\frac{(1+\theta-\gamma)}{p-r}\leq \frac{N}{N-r},$
so that
$$r\leq \frac{N[p-1-\theta + \gamma]}{N-1-\theta+\gamma}.$$
Then $p>p_{0}$ implies
$$\frac{N[p-1-\theta + \gamma]}{N-1-\theta+\gamma}>1.$$
We obtain
$$\displaystyle\int_{\Omega} u_{n}^{r^{\ast}}dx \leq
c	\left(\displaystyle\int_{\Omega}  (1+u_{n})^{r^{\ast}}dx \right)^{\frac{(p-r)r^{\ast}}{p-r}}\leq c+c \left(\displaystyle\int_{\Omega} u_{n}^{r^{\ast}}dx \right)^{\frac{(p-r)r^{\ast}}{pr}}.$$
From the above inequality, by Young's inequality with $\varepsilon$, we see that
$$ 	\int_{\Omega}  | u_{n}|^{r^{\ast}}dx \leq c .$$
Which together with \eqref{37} and $\frac{r}{p-r}(\theta+1-\gamma)\leq r^{\ast}$
implies
$$	\int_{\Omega}  |\nabla u_{n}|^{r} dx\leq c \,\,\,\,\,, r< \frac{N[p-1-\theta + \gamma]}{N-1-\theta+\gamma}.$$
Just as in the proof of part I in the Theorem \ref{t3}, we can conclude that $u$ is a distributional solution of \eqref{S1}.\\
\textbf{Part 2.} Let
$$\frac{N(1-\gamma)+\gamma}{N-1}<s<\frac{N(p-\theta-1)+p\gamma}{N-p}.$$
We can show that $u$ is an entropy solution \eqref{S1} by the same method in part \textbf{b} of Theorem \ref{t3}. \\
\textbf{Part 3.} Let $s\geq \frac{N(p-\theta-1)+p\gamma}{N-p}$. Since  $p>p_{0},$ this implies
$$\frac{N(p-\theta-1)+p\gamma}{N-p}>\frac{1+\theta -\gamma}{p-1}.$$
Therefore, the proof of this part is the same as the proof of part \textbf{a} in Theorem$\ref{t3}$.

\end{document}